\documentclass[12pt,thmsa]{article}
\usepackage{amsfonts}
\usepackage{amssymb}
\usepackage{graphicx}
\usepackage{amssymb,amsmath}
\newtheorem{teo}{Theorem}[section]

\newtheorem{obs2}[teo]{Remark}

\newtheorem{tea}{Theorem}[subsection]

\newtheorem{no2}[teo]{Note}

\newtheorem{no3}[tea]{Note}

\newcommand{\Q}{{\mathbb{Q}}}

\def\timehm{\count31=\time \count32=\count31 \divide\count31 by 60
\number\count31 \multiply\count31 by 60 \advance\count32 by
-\count31 :\ifnum\count32<10 0\fi \number\count32}

\def\ideal#1{<\kern-2pt #1\kern-2pt >}

\begin{document}
\title{{\bf  Demostraci\'{o}n simplificada de la Conjetura de Modularidad de Serre (carta)}
}
\author{Luis Dieulefait
%\thanks{Research partially supported by
 % MTM2012-33830 and by an ICREA Academia Research Prize}
\\
Universitat de Barcelona\\
G.V. de les Corts Catalanes 585\\
08007 - Barcelona, Spain \\
e-mail: ldieulefait@ub.edu\\
%Date:  (Preliminary Version)
 }
\date{\empty}

\maketitle

\vskip -20mm

\begin{abstract} We explain in this letter how using a recent Modularity Lifting Theorem proved by Lue Pan the proofs of Serre's Modularity Conjecture over $\Q$ given by Khare-Wintenberger and the author can be greatly simplified. The main difference with the previous proofs is that neither the process of ``weight reduction" developed by Khare in his proof of the level $1$ case nor the alternative method of ``weight reduction via Galois conjugation" developed by the author are required: weight reduction can now be obtained just with two applications of Pan's result, together with a Lemma that guarantees  that the residual representations in the last two steps are not in the bad-dihedral case.
\end{abstract}

\newpage

\begin{quote}
\it  A Jean-Pierre Wintenberger (1954-2019), in Memoriam.
\end{quote}

\section{Carta datada en Septiembre de  2018}

Hola Basil, \\

Te hago notar que, usando el teorema de L. Pan, en la versi\'{o}n que \'{e}l dijo que
pronto demostrar\'{a}, en la cual quita la condici\'{o}n tediosa local en $p$ (sobre el
cociente de los dos caracteres en la representaci\'{o}n residual), al menos para $p > 3$, se puede probar la conjetura de Serre sobre $\Q$ de manera muy simple. De
hecho, la fase de ``killing the level" es muy sencilla y ya est\'{a} explicada en mi
paper ``Remarks on Serre's modularity conjecture", y una vez hecho esto, se
acaba con el teorema de Pan en s\'{o}lo dos pasos m\'{a}s! \\
Te recuerdo la parte de ``killing the level": como el teorema de Kisin ya ha sido
mejorado por otros autores (Y. Hu, F. Tan y S. Tung), quitando precisamente la condici\'{o}n tediosa sobre
el cociente de los dos caracteres de la representaci\'{o}n residual local en $p$, tal
como anticip\'{e} que pasar\'{i}a en la secci\'{o}n final de mi paper (secci\'{o}n 7),  ``killing
ramification" se hace de la forma m\'{a}s simple imaginable: se introduce un Good
Dihedral prime $q$ en el nivel para garantizar ``large residual images" y se
matan uno tras otro los primos $p_i$ del nivel, es decir, se reduce m\'{o}dulo $p_i$ y se
toma un levantamiento minimal sin $p_i$ en el nivel. Hecho esto, la conjetura de
Serre queda reducida al caso de nivel $q^2$, peso $k \geq 2$. 
Es decir, se aplica el m\'{e}todo tal como fuera explicado en mi paper, pero del modo simplificado que se mencionaba en la Secci\'{o}n 7 del mismo gracias a que el teorema de Kisin ya ha sido mejorado.\\

Observaci\'{o}n: en realidad, estoy suponiendo que el nivel es impar. Si fuera par,
para matar el primo $2$ del nivel, se procede como en los papers de Khare-Wintenberger, la idea es la siguiente: despu\'{e}s de
haber hecho lo anterior (introducir el Good Dihedral prime $q$ en el nivel y
matar todos los primos impares del nivel) tomando un levantamiento de peso
$2$ modulo un primo auxiliar $r$ se puede suponer que se est\'{a} en un caso de peso
$2$ (el primo $r$ ser\'{a} un primo que est\'{a} ahora en el nivel pero que se mata luego
de la forma sencilla ya indicada en el p\'{a}rrafo anterior, basta con cogerlo menor que una cota apropiada para que el good dihedral prime $q$ garantice imagen residual grande modulo $r$), y una vez hecho esto, el
\'{u}nico problema t\'{e}cnico para matar el $2$ del nivel es que los ``$2$-adic modularity
lifting theorems" (el de Kisin es el que ellos aplican) requieren que sea una
representaci\'{o}n potentially Barsotti-Tate (no son v\'{a}lidos en el caso
semiestable no-cristalino). Para resolver este problema t\'{e}cnico, Khare-Wintenberger hacen
un truco que podr\'{i}amos llamar ``three-two trick": Si est\'{a}n en el caso malo, o
sea con un sistema compatible que localmente en el $2$ est\'{a} en el caso
semistable (up to twist), prueban que se puede hacer una congruencia
modulo $3$ con otro sistema cuyo tipo local en el $2$ sea tal que la representaci\'{o}n
$2$-adica del nuevo sistema es potentially Barsotti-Tate. Hecho esto, ya se puede
usar el teorema $2$-adico de Kisin para matar el $2$ del nivel, y se tiene una
congruencia con una representaci\'{o}n de conductor impar y peso $2$ o $4$. Es en este momento d\'{o}nde podemos matar el primo $r$ del nivel tal como hicimos con el resto de primos impares del nivel. \\

Una vez que se acaba la fase de killing the level, que es muy sencilla, ya
estamos en el caso de nivel $q^2$ y peso $k > 1$. Aqu\'{i}, si suponemos la mejora del teorema
de Pan (para $p > 3$) que \'{e}l anunci\'{o} en su charla del congreso de Rio de Janeiro
(tal como dije al principio: un teorema que no requiera la condici\'{o}n tediosa local en $p$ para $p > 3$, con
lo que queda un teorema como el de Skinner-Wiles, es decir un teorema de
modularidad para el caso residualmente reducible, pero que cubre tambi\'{e}n el caso no-ordinario, y que tiene como \'{u}nicas condiciones las que son necesarias para que
la representaci\'{o}n $p$-adica considerada pueda ser modular, o sea las
condiciones de la conjetura de Fontaine-Mazur-Langlands) se puede acabar en
dos pasos: reducir modulo $q$ y, si fuera residualmente irreducible y con imagen large, tomar un levantamiento minimal para as\'{i} pasar a una congruencia con un sistema de
nivel $1$ y peso mayor que $1$ (aqu\'{i} hay que tener en cuenta que al reducir modulo $q$
estamos perdiendo el Good Dihedral prime y ya no tenemos control sobre la
imagen residual): si en este paso la imagen residual es ``large" se aplica de
nuevo el Teorema de Kisin (como en mi paper ya citado) para garantizar que la modularidad se preserva en esta congruencia (se propaga de la de nivel $1$ a la de nivel $q^2$). Si la representaci\'{o}n
residual es reducible, se aplica el teorema de L. Pan (en la versi\'{o}n mejorada) y se concluye modularidad. Queda
un tercer caso, que es el caso donde la imagen residual es bad-dihedral (o
sea, irreducible pero tal que su restricci\'{o}n al cuerpo cuadr\'{a}tico ramificado s\'{o}lo
en $q$ es reducible), pero en un lema al final de este mensaje probar\'{e}
que este tercer caso no puede ocurrir. Notese que queremos eliminar este caso para poder aplicar el Teorema de Kisin y para poder aplicar el resultado de existencia de un levantamiento minimal.\\
En el paso siguiente, que s\'{o}lo se lleva a cabo si la representaci\'{o}n modulo $q$ es
irreducible, se coge el levantamiento minimal de nivel 1 y peso $k  > 1$, y el
sistema compatible que lo contiene, y se reduce modulo $5$. Tras un apropiado
twist, podemos suponer que el peso de Serre de esta representaci\'{o}n m\'{o}dulo $5$
es $k =2, 4$ o $6$. En los tres casos, la conjetura de Serre la damos por cierta: estos
tres casos base se reducen f\'{a}cilmente (usando la maquinaria de modularity
lifting theorems, minimal lifts, compatible systems...) a resultados de Serre (y
Tate) m\'{o}dulo $p=3$ para los casos  $k = 2 $ y $4$ y a resultados de Schoof para $k =6$
(sobre variedades abelianas semiestables con buena reducci\'{o}n fuera del $5$). De
hecho, son tres casos de la conjetura de Serre que fueron probados en los
primeros papers sobre conjetura de Serre en trabajos de Khare-Wintenberger y
m\'{i}os de 2003 y 2004. As\'{i} que los consideramos como casos ya probados. Como
no hay formas modulares cuspidales con nivel $1$ y peso $2, 4$ o $6$, esto significa que la
representaci\'{o}n m\'{o}dulo $5$ tiene que ser reducible 
(de hecho, aplicando el lema al final de esta carta se ve en particular que no puede ser bad-dihedral), con lo cual se concluye
modularidad del sistema compatible aplicando una vez m\'{a}s el teorema de L.
Pan. Y esto acaba la demostraci\'{o}n de la Conjetura de Serre. \\

S\'{o}lo falta demostrar el lema que explica porqu\'{e} al reducir m\'{o}dulo $q$ la
representaci\'{o}n residual no puede ser bad-dihedral: \\
Demostraci\'{o}n del lema: Si la representaci\'{o}n residual es bad-dihedral, tendr\'{a} peso de
Serre $k > 1$ (twisteando si fuera necesario, podemos suponer  $ 1 < k \leq  q + 1 $), nivel
de Serre $1$, y la caracter\'{i}stica del cuerpo es el primo $q$. Utilizando
resultados de Ribet en los que este caso de imagen bad-dihedral se estudia en
detalle, es bien sabido, y se aplic\'{o} en reiteradas ocasiones en las demostraciones
de la conjetura de Serre de la pasada d\'{e}cada (tanto en la de Khare-Wintenberger como en la m\'{i}a, y de hecho tambi\'{e}n se aplic\'{o} en el resultado
anterior de Dieulefait-Manoharmayum sobre modularidad de variedades de
Calabi-Yau r\'{i}gidas) que para que la representaci\'{o}n sea bad-dihedral se tiene
que tener la siguiente relaci\'{o}n entre el primo $q$ y el peso de Serre $k$:
o bien $q = 2 k - 3$, o bien $q = 2 k - 1$. \\
La primera se corresponde con el caso en que la inercia en $q$ act\'{u}a a trav\'{e}s de
caracteres fundamentales de niveau $2$ y la segunda  con el caso de
niveau $1$, es decir, potencias del caracter ciclot\'{o}mico modulo $q$. 
Ahora bien, en el caso en que estamos, en el cual el nivel de Serre es $1$, un
an\'{a}lisis m\'{a}s fino  da lugar al ``Lema de Wintenberger" (cuya
demostraci\'{o}n aparece en el \'{i}tem $1$ del lema 6.2 del paper de C. Khare ``Serre's
modularity conjecture: the level 1 case" publicado en el Duke Math. Journal) en el que
se demuestra que el caso bad-dihedral con $q = 2 k -3 $ (es decir, el caso de
caracteres fundamentales de niveau $2$) no puede ocurrir. Con lo cual, la \'{u}nica
posibilidad es $ q = 2 k -1$. Pero observemos que como el nivel de Serre es 1,
nuestra representación residual (debido a que es una representaci\'{o}n impar) tiene que tener un
valor de $k$ par, con lo cual la igualdad $q = 2 k - 1$ implica que $q$ es congruente
con $3$ modulo $4$. Pero esto contradice la definici\'{o}n misma de Good Dihedral
prime, puesto que parte de la definici\'{o}n impone que $q$ ha de ser congruente
con $1$ m\'{o}dulo $8$ (vease el paper ``Serre's modularity conjecture (I)", de Khare-
Wintenberger). Por lo tanto el lema queda probado, y es v\'{a}lido en general
para cualquier primo $p$ congruente con $1$ modulo $4$: si el nivel de Serre es $1$
en una representaci\'{o}n modulo $p$ para un tal primo $p$, entonces la
representaci\'{o}n no puede ser bad-dihedral (por eso el comentario que hice antes de que este lema tambi\'{e}n podr\'{i}a aplicarse al reducir modulo $5$ para deducir que la representaci\'{o}n residual no es bad-dihedral, puesto que $5$ es congruente con $1$ modulo $4$). \\

Bueno, habiendo acabado por completo esta demostraci\'{o}n breve de la
conjetura de Serre, queda s\'{o}lo un comentario final: 
 Cu\'{a}l es la principal ventaja de esta demostraci\'{o}n m\'{a}s corta? Que es mucho
m\'{a}s sencilla. El mecanismo de inducci\'{o}n doble (en el nivel y en el peso) que aparece en la demostraci\'{o}n dada por Khare-Wintenberger con la casu\'{i}stica que incluye es claramente m\'{a}s complicado. Por otro
lado, en mi paper ``Remarks on Serre's modularity conjecture" ya se utiliza el resultado
m\'{a}s fuerte de Kisin para hacer el proceso de killing the level m\'{a}s sencillo,
pero una vez reducidos a casos de nivel $1$, hay que recordar que para reducir el
peso de Serre tanto el m\'{e}todo de reducci\'{o}n de peso de Khare (en su paper publicado en
el Duke Math. Journal sobre el caso de nivel $1$) como mi m\'{e}todo alternativo en mi paper ``The level 1 case of Serre's
conjecture revisited", publicado en Rendiconti dei Lincei (la ventaja de este
paper era que evitaba el uso de primos auxiliares que dividieran a $p-1$, cosa
que hac\'{i}a Khare en su paper, mediante el truco de ``weight reduction via Galois
conjugation") utilizaban versiones refinadas del postulado de Bertrand, con lo
cual se trataba en ambos casos de un m\'{e}todo inductivo sobre los primos
utilizando t\'{e}cnicas anal\'{i}ticas. La demostraci\'{o}n en esta carta, gracias al
teorema de L. Pan, ejecuta la reducci\'{o}n de peso en un s\'{o}lo paso, simplemente
movi\'{e}ndonos hasta el primo $p=5$. \\

Post Scriptum(January 2019): As expected, L. Pan proved his result in the stronger version required for the argument in this letter to work, this is the main result in his paper: ``The Fontaine-Mazur conjecture in the residually reducible case" (L. Pan, available at ArXiv). Therefore, the short proof of Serre`s conjecture described in this letter is now unconditional.

\end{document}